% To: shlhetal@math.huji.ac.il
% Subject: none
% Date: Fri, 30 Apr 2004 10:28:23 +0200
% From: Tomas Jech <jech@math.cas.cz>
% X-sliced-and-diced-by: 'savemail' 1.3b, Feb 2003

%
\magnification=\magstep1
\input amstex
\documentstyle{amsppt}

% The next line affects recent versions of amsppt.sty,
% making them emulate the pre-1997 spacing.

\ifx\undefined\headlineheight\else\pageheight{47pc}\vcorrection{-1pc}\fi
 
\hsize=30truecc
\NoBlackBoxes
\leftheadtext{Thomas Jech and Saharon Shelah}
\def\pf{\hfill $\qed$}
\def\c{\cite}

\def\ifof{\text{if and only if }}
\def\seq{\text{Seq }}
\def\doc{\dot{C}}
\def\dof{\dot{f}}
\def\dox{\dot{X}}

\topmatter
\title
Simple Complete Boolean Algebras
\endtitle
\author
Thomas Jech$^1$ and Saharon Shelah$^{2,3}$
\endauthor
\affil 
The Pennsylvania State University  \\
The Hebrew University of Jerusalem and Rutgers University
\endaffil
\thanks
$^{1,2}$Supported in part by the National Science Foundation grant 
DMS--98-02783 and DMS--97-04477.
\endthanks
\address
Department of Mathematics, The Pennsylvania State University,  \newline
218 McAllister Bldg., University Park, PA 16802, U.S.A. \newline
Institute of Mathematics, The Hebrew University of Jerusalem,  \newline
91904 Jerusalem, Israel.
\endaddress
\email jech\@math.psu.edu,
	 shelah\@math.rutgers.edu 
\endemail
\abstract
For every regular cardinal $\kappa$ there exists a simple 
complete Boolean algebra with $\kappa$ generators.
\endabstract
\subjclass 
\endsubjclass
\endtopmatter

\document
\baselineskip=16truept
\footnote"{}"{$^3$Paper number 694.}

\subhead
1. \ Introduction
\endsubhead

A complete Boolean algebra is {\it simple} if it is atomless and 
has no nontrivial proper atomless complete subalgebra.  The problem
of the existence of simple complete Boolean algebras was first
discussed in 1971 by McAloon in \c8.  Previously, in \c7,
McAloon constructed a {\it rigid} complete Boolean algebra;
it is easily seen that a simple complete Boolean algebra is 
rigid.  In fact, it has no non-trivial one-to-one complete
endomorphism \c1.  Also, if an atomless complete algebra is
not simple, then it contains a non-rigid atomless complete 
subalgebra \c2.

McAloon proved in \c8 that an atomless complete algebra $B$
is simple \ifof it is rigid and {\it minimal}, i.e.
the generic extension by $B$ is a minimal extension of the
ground model.  Since Jensen's construction \c5 yields a 
definable real of minimal degree over $L$, it shows that a
simple complete Boolean algebra exists under the assumption
$V = L$.  McAloon then asked whether a rigid minimal algebra
can be constructed without such assumption.

In \c{10}, Shelah proved the existence of a rigid complete Boolean
algebra of cardinality $\kappa$ for each regular cardinal $\kappa$
such that $\kappa^{\aleph_0} = \kappa$.  Neither McAloon's nor
Shelah's construction gives a minimal algebra.

In \c9, Sacks introduced perfect set forcing, to produce a
real of minimal degree.  The corresponding complete Boolean 
algebra is minimal, and has $\aleph_0$ generators.  Sacks' 
forcing generalizes to regular uncountable cardinals $\kappa$
(cf.\c6), thus giving a minimal complete Boolean algebra with
$\kappa$ generators.  The algebras are not rigid however.

Under the assumption $V = L$, Jech constructed in \c3 a simple
complete Boolean algebra of cardinality $\kappa$, for every
regular uncountable cardinal that is not weakly compact (if 
$\kappa$ is weakly compact, or if $\kappa$ is singular and GCH 
holds, then a simple complete Boolean algebra does not exist).

In \c4, we proved the existence of a simple complete Boolean algebra
(in ZFC).  The algebra is obtained by a modification of Sacks'
forcing, and has $\aleph_0$ generators (the forcing produces a
definable minimal real).  The present paper gives a construction of a
simple complete Boolean algebra with $\kappa$ generators, for every
regular uncountable cardinal $\kappa$.

\proclaim
{Main Theorem}  
Let $\kappa$ be a regular uncountable cardinal.
There exists a forcing notion $P$ such that the complete Boolean
algebra $B = B(P)$ is rigid, $P$ adds a subset of $\kappa$ 
without adding any bounded subsets, and for every $X \in V[G]$
(the $P$-generic extension), either $X \in V$ or $G \in V[X]$.
Consequently, $B$ is a simple complete Boolean algebra with $\kappa$
generators.
\endproclaim

The forcing $P$ is a modification of the generalization of Sacks'
forcing described in \c{6}.

\subhead
2. \ Forcing with perfect $\kappa$-trees
\endsubhead

For the duration of the paper let $\kappa$ denote a regular
uncountable cardinal, and set Seq $=\bigcup _{\alpha < \kappa}
\,^\alpha 2$. 

\definition
{Definition 2.1}  (a) \ If $p \subseteq \seq$ and $s \in p$, say 
that $s$ {\it splits} in $p$ if $s^\frown 0 \in p$ and $s^\frown 1 
\in p$.

(b) \ Say that $p \subseteq \seq$ is a {\it perfect tree} if:
\roster
\item"{(i)}"  If $s \in p$, then $s{\upharpoonright} \alpha \in p$
	for every $\alpha$.
\item"{(ii)}" If $\alpha < \kappa$ is a limit ordinal, $s \in \,^\alpha 2$,
	and $s{\upharpoonright} \beta \in p$ for every $\beta < \alpha$,
	then $s \in p$.
\item"{(iii)}" If $s \in p$, then there is a $t \in p$ with $t 
	\supseteq s$ such that $t$ splits in $p$.
\endroster
Our definition of perfect trees follows closely \c6, with one
exception: unlike \c6, Definition 1.1.(b)(iv), the splitting
nodes of $p$ need not be closed.
\enddefinition

We consider a notion of forcing $P$ that consists of (some)
perfect trees, with the ordering $p \leq q$ iff $p \subseteq q$.
Below we formulate several properties of $P$ that guarantee that
the proof of minimality for Sacks forcing generalizes to forcing
with $P$.

\definition
{Definition 2.2} (a) \ If $p$ is a perfect tree and $s \in p$, set
$p_s = \{t \in p : s \subseteq t \text{ or } t \subseteq s\}$; 
$p_s$ is a {\it restriction} of $p$.  A set $P$ of perfect trees
is {\it closed under restrictions} if for every $p \in P$ and
every $s \in p$, $p_s \in P$.  If $p_s = p$, then $s$ is a {\it stem}
of $p$.

(b) \ For each $s \in \seq$, let $o(s)$ denote the domain (length)
of $s$.  If $s \in p$ and $o(s)$ is a successor ordinal, $s$ is
a {\it successor node} of $p$; if $o(s)$ is a limit ordinal,
$s$ is a {\it limit node} of $p$.  If $s$ is a limit node of
$p$ and $\{ \alpha < o(s) : s \upharpoonright \alpha \text{ 
splits in }p\}$ is cofinal in $o(s)$, $s$ is a {\it limit of
splitting nodes}.

(c) \ Let $p$ be a perfect tree and let $A$ be a nonempty set 
of mutually incomparable successor nodes of $p$.  If for each
$s \in A$, $q(s)$ is a perfect tree with stem $s$ and $q(s) \leq
p_s$, let
$$
	q = \{ t \in p : \text{ if } t \supseteq s \text{ for 
	some } s \in A \text{ then } t \in q(s)\}
$$
We call the perfect tree $q$ the {\it amalgamation} of
$\{q(s) : s \in A\}$ into $p$.  A set $P$ of perfect trees is
{\it closed under amalgamations} if for every $p \in P$, every
set $A$ of incomparable successor nodes of $p$ and every
$\{q(s) : s \in A\} \subset P$ with $q(s) \leq p_s$, the
amalgamation is in $P$.
\enddefinition

\definition
{Definition 2.3} (a) \ A set $P$ of perfect trees is $\kappa$-{\it
closed} if for every $\gamma < \kappa$ and every decreasing sequence
$	\langle p_{\alpha} : \alpha < \gamma\rangle \text{ in }
	P,\, \bigcap_{\alpha < \gamma} p_{\alpha} \in P.
$

(b) \ If $\langle p_{\alpha} : \alpha < \kappa\rangle$ is a 
decreasing sequence of perfect trees such that
\roster
\item"{(i)}"  if $\delta$ is a limit ordinal, then $p_{\delta}
	= \bigcap_{\alpha < \delta} p_{\alpha}$, and
\item"{(ii)}" for every $\alpha$, $p_{\alpha + 1} \cap \,^\alpha 2
	= p_{\alpha} \cap \,^{\alpha}2$, 
\endroster
then $\langle p_{\alpha}: \alpha < \kappa\rangle$ is called a 
{\it fusion sequence}.  A set $P$ is {\it closed under fusion}
if for every fusion sequence
$	\langle p_{\alpha} : \alpha < \kappa\rangle \text{ in }
	P,\, \bigcap_{\alpha < \kappa} p_{\alpha} \in P\,.
$
\enddefinition

The following theorem is a generalization of Sacks' Theorem from
\c{9} to the uncountable case:

\proclaim
{Theorem 2.4}  Let $P$ be a set of perfect trees and assume that
$P$ is closed under restrictions and amalgamations, $\kappa$-closed,
and closed under fusion.  If $G$ is $P$-generic over $V$, then $G$
is minimal over $V$; namely if $X \in V[G]$ and $X \notin V$, 
then $G \in V[X]$.  Moreover, $V[G]$ has no new bounded subsets
of $\kappa$, and $G$ can be coded by a subset of $\kappa$.
\endproclaim

\demo
{Proof}  The proof follows as much as in \c9.  Given a name
$\dox$ for a set of ordinals and a condition $p \in P$ that 
forces $\dox \notin V$, one finds a condition $q \leq p$ and a
set of ordinals $\{\gamma_s : s$ splits in $q\}$ such that
$q_{s^{\frown}0}$ and $q_{s^{\frown}1}$  both decide $\gamma_s
\in \dox$, but in opposite ways.  The generic branch can then
be recovered from the interpretation of $\dox$.

To construct $q$ and $\{\gamma_s\}$ one builds a fusion sequence
$\{p_{\alpha} : \alpha < \kappa\}$ as follows.  Given $p_{\alpha}$,
let $Z = \{s \in p_{\alpha} : o(s) = \alpha$ and $s$ splits in
$p_{\alpha}\}$.  For each $s \in Z$, let $\gamma_s$ be an ordinal
such that $(p_{\alpha})_s$ does not decide $\gamma_s \in \dox$.
Let $q(s^{\frown}0) \leq (p_{\alpha})_{s^{\frown}0}$ and
$q(s^{\frown}1) \leq (p_{\alpha})_{s^{\frown}1}$ be conditions
that decide $\gamma_s \in \dox$ in opposite ways.  Then let
$p_{\alpha + 1}$ be the amalgamation of $\{q(s^{\frown}i) :
s \in Z$ and $i = 0,1\}$ into $p_{\alpha}$.  Finally, let
$q = \bigcap_{\alpha < \kappa} p_{\alpha}$.  \hfill$\qed$
\enddemo

In \c6 it is postulated that the splitting nodes along any branch
of a perfect tree form a closed unbounded set.  This guarantees
that the set of all such trees is $\kappa$-closed and closed under
fusion (Lemmas 1.2 and 1.4 in \c6).  It turns out that a less
restrictive requirement suffices.

\definition
{Definition 2.5}  Let $S \subset \kappa$ be a stationary set.  A
perfect tree $p \in P$ is $S$-{\it perfect} if whenever $s$ is a
limit of splitting nodes of $p$ such that $o(s) \in S$, then
$s$ splits in $p$.
\enddefinition

\proclaim
{Lemma 2.6}  (a) \ If $\langle p_{\alpha} : \alpha < \gamma\rangle$,
$\gamma < \kappa$, is a decreasing sequence of $S$-perfect trees,
then $\bigcap_{\alpha < \gamma} p_{\alpha}$ is a perfect
tree.

(b) \ If $\langle p_{\alpha} : \alpha < \kappa\rangle$ is a fusion
sequence of $S$-perfect trees, then $\bigcap_{\alpha < \kappa}
p_{\alpha}$ is a perfect tree. 
\endproclaim

\demo
{Proof}  (a) \ Let $p = \bigcap_{\alpha < \gamma} p_{\alpha}$.
The only condition in Definition 2.1 (b) that needs to be verified
is (iii): for every $s\in p$ find $t \supseteq s$ that splits in $p$.
First it is straightforward to find a branch $f \in \,^\kappa 2$
through $p$ such that $s$ is an initial segment of $f$.

Second, it is equally straightforward to see that for each $\alpha
< \gamma$, the set of all $\beta$ such that $f \upharpoonright
\beta$ splits in $p_{\alpha}$ is unbounded in $\kappa$.  Thus for
each $\alpha < \gamma$ let $C_{\alpha}$ be the closed unbounded
set of all $\delta$ such that $f \upharpoonright \delta$ is a
limit of splitting nodes in $p_{\alpha}$.  Let $\delta \geq o(s)$
be an ordinal in $\bigcap_{\alpha < \gamma} C_{\alpha} 
\cap S$.  Then for each $\alpha < \gamma$, $t = f\upharpoonright
\delta$ is a limit of splitting nodes in $p_{\alpha}$, and hence $t$
splits in $p_{\alpha}$.  Therefore $t$ splits in $p$.

(b) \ Let $p = \bigcap_{\alpha < \kappa} p_{\alpha}$
and again, check (iii).  Let $s \in p$, and let $f \in \,^\kappa 2$
be a branch trough $p$.  For each $\alpha < \gamma$ let $C_{\alpha}$
be the club of all $\delta$ such that $f \upharpoonright\delta$ 
is a limit of splitting nodes in $p_{\alpha}$.  Let $\delta \geq
o(s)$ be an ordinal in $\Delta_{\alpha<\kappa} C_{\alpha} \cap S$ and let
$t = f \upharpoonright \delta$.  If $\alpha < \delta$, then
$t$ splits in $p_{\alpha}$, and therefore $t$ splits in $p_{\delta}$.
Since $p_{\delta + 1}\cap \,^\delta 2 = p_{\delta} \cap \,^\delta 2$,
we have $t \in p_{\delta + 1}$, and since $p_{\delta + 1}$ is
$S$-perfect, $t$ splits in $p_{\delta + 1}$.  If $\alpha >
\delta + 1$, then $p_{\alpha} \cap \,^{\delta + 1}2 = p_{\delta+1}
\cap \,^{\delta + 1}2$, and so $t$ splits in $p_{\alpha}$.  Hence
$t$ splits in $p$.  \pf
\enddemo

This is trivial, but note that the limit condition $p$ (in both
(a) and (b)) is not only perfect but $S$-perfect as well.

\subhead
3. \ The notion of forcing for which $B(P)$ is rigid
\endsubhead

We now define a set $P$ of perfect $\kappa$-trees that is
closed under restrictions and amalgamations, $\kappa$-closed,
and closed under fusion, with the additional property that 
the complete Boolean algebra $B(P)$ is rigid.  That completes
a proof of Main Theorem.

Let $S$ and $S_{\xi}$, $\xi < \kappa$, be mutually disjoint
stationary subsets of $\kappa$, such that for all $\xi < \kappa$,
if $\delta \in S_{\xi}$, then $\delta > \xi$.

\definition
{Definition 3.1}  The forcing notion $P$ is the set of all
$p \subseteq \seq$ such that
\roster
\item $p$ is a perfect tree;
\item $p$ is $S$-perfect, i.e. if $s$ is a limit of splitting nodes
	of $p$ and $o(s)\in S$, then $s$ splits in $p$;
\item For every $\xi < \kappa$, if $s$ is a limit of splitting
	nodes of $p$ with $o(s) \in S_{\xi}$ and if $s(\xi) = 0$
	then $s$ splits in $p$.
\endroster

The set $P$ is ordered by $p \leq q$ iff $p \subseteq q$.
\enddefinition

Clearly, $P$ is closed under restrictions and amalgamations.  By
Lemma 2.6, the intersection of either a decreasing short sequence
or of a fusion sequence in $P$ is a perfect tree, and since both
properties (2) and (3) are preserved under arbitrary intersections,
we conclude that $P$ is also $\kappa$-closed and closed under
fusion.

We conclude the proof by showing that $B(P)$ is rigid.

\proclaim
{Lemma 3.2}  If $\pi$ is a nontrivial automorphism of $B(P)$, then
there exist conditions $p$ and $q$ with incomparable stems such
that $\pi(p)$ and $q$ are compatible (in $B(P)$).
\endproclaim

\demo
{Proof} Let $\pi$ be a nontrivial automorphism.  It is easy to find a
nonzero element $u \in B$ such that $\pi(u) \cdot u = 0$.  Let $p_1
\in P$ be such that $p_1 \leq u$, and let $q_1 \in P$ be such that
$q_1 \leq \pi(p_1)$.  As $p_1$ and $q_1$ are incompatible, there
exists some $t \in q_1$ such that $t \notin p_1$.  Let $q = (q_1)_t$.
Then let $p_2 \in P$ be such that $p_2\leq \pi^{-1}(q)$, and again,
there exists some $s \in p_2$ such that $s \notin q$.  Let $p =
(p_2)_s$.  Now $s$ and $t$ are incomparable stems of $p$ and $q$, and
$\pi(p) \leq q$.  \pf
\enddemo

To prove that $B(P)$ has no nontrivial automorphism, we introduce
the following property $\varphi(\xi)$.

\definition
{Definition 3.3}  Let $\xi < \kappa$; we say that $\xi$ has
property $\varphi$ if and only if for every function $f : \kappa
\to 2$ there exist a function $F : \seq \to 2$ in $V$ and a club
$C \subset \kappa$ such that for every $\delta \in C \cap 
S_{\xi}$, $f(\delta) = F(f \upharpoonright \delta)$.
\enddefinition

\proclaim
{Lemma 3.4}  Let $t_0 \in \seq$ and let $\xi = o(t_0)$.
\roster
\item"{(a)}"  Every condition with stem $t_0^{\frown} 0$ forces
	$\neg\varphi(\xi)$.
\item"{(b)}"  Every condition with stem $t_0^{\frown} 1$ forces
	$\varphi(\xi)$.
\endroster
\endproclaim

\demo
{Proof}  (a) \ Let $\dof$ be the name for the generic branch
$f_G : \kappa \to 2$ (i.e. $f_G = \bigcup \{s \in \seq : s \in p
\text{ for all } p \in G\}$); this will be the counterexample for
$\varphi(\xi)$.  Let $F$ be a function, $F : \seq \to 2$, let 
$\doc$ be a name for a club and let $p \in P$ be such that 
$t_0^{\frown}0$ is a stem of $p$.  We shall find a $\delta \in
S_{\xi}$ and $q \leq p$ such that $q \Vdash (\delta\in \doc \text{ and }
\dof(\delta) \neq F (\dof \upharpoonright \delta))$.

We construct a fusion sequence $\langle p_{\alpha} : \alpha < \kappa\rangle$,
starting with $p$, so that for each $\alpha$, if $s \in p_{\alpha +
1}$ and $o(s) = \alpha + 1$, then $(p_{\alpha + 1})_s$ decides the
value of the $\alpha$th element of $\doc$; we call this value
$\gamma_s$.  (We obtain $p_{\alpha + 1}$ by amalgamation into
$p_{\alpha}$.)  Let $r = \bigcap_{\alpha < \kappa}
p_{\alpha}$.

Let $b$ be a branch through $r$, and let $s_{\alpha} = b \upharpoonright
\alpha$ for all $\alpha$.  There exists a $\delta \in S$ such that $s_{\delta}$
is a limit of splitting nodes of $r$, and such that for every $\alpha <
\delta$, $\gamma_{s_{\alpha + 1}} < \delta$.  Since $s_{\delta}(\xi)
= 0$, $s_{\delta}$ splits in $r$, and $r_{s_{\delta}} \Vdash \delta
\in \doc$.

Now if $F(s_{\delta}) = i$, it is clear that $g = r_{s_{\delta}^{\frown}
(1 - i)}$ forces $\dof \upharpoonright \delta = s_{\delta}$ and
$\dof(\delta) = 1 - i$.

(b) \ Let $\dof$ be a name for a function from $\kappa$ to $2$,
and let $p$ be a condition with stem $t_0^{\frown} 1$ that forces
$\dof \notin V$ ($\varphi(\xi)$ holds trivially for those $f$
that are in $V$).  We shall construct a condition $q \leq p$
and collections $\{h_s : s \in Z\}$ and $\{i_s : s \in Z'\}$,
where $Z$ is the set of all limits of splitting nodes in $q$
and $Z' = \{s \in Z : o(s) \in S_{\xi}\}$, such that
\enddemo

{\bf (3.5)} 
\roster
\item"{(i)}"     For each $s \in Z$, $h_s \in \seq$ and $o(h_s)
		= o(s)$; if $o(s) = \alpha$, then $q_s\Vdash \dof
		\upharpoonright \alpha = h_s$.
\item"{(ii)}"    If $s,t \in Z$, $o(s) = o(t) = \alpha$, and $s \neq t$,
		then $h_s \neq h_t$.
\item"{(iii)}"   For each $s \in Z'$, $i_s = 0$ or $i_s = 1$; if
		$o(s) = \delta$, then $q_s \Vdash \dof(\delta) = i_s$.
\endroster

Then we define $F$ by setting $F(h_s) = i_s$, for all $s \in Z'$
(and $F(h)$ arbitrary for all other $h \in \seq$); this is possible
because of (ii).  We claim that $q$ forces that for some club $C$,
$\dof(\delta) = F(\dof \upharpoonright \delta)$ for all $\delta
\in C \cap S_{\xi}$.  (This will complete the proof.)  

To prove the claim, let $G$ be a generic filter with $q \in G$,
let $g$ be the generic branch $(g = \bigcup \{s : s \in p \text{ for 
all } p \in G\})$, and let $f$ be the $G$-interpretation of $\dof$.
Let $C$ be the set of all $\alpha$ such that $g\upharpoonright \alpha$
is the limit of splitting nodes in $q$.  If $\delta \in C \cap 
S_{\xi}$, let $s = g \upharpoonright \delta$; then $s \in Z'$,
$f \upharpoonright\delta = h_s$ and $f(\delta) = i_s$.  It 
follows that $f(\delta) = F(f \upharpoonright \delta)$.

To construct $q$, $h_s$ and $i_s$, we build a fusion sequence
$\langle p_{\alpha} : \alpha < \kappa\rangle$ starting with $p_0$.  We
take $p_{\alpha} = \bigcap_{\beta <\alpha} p_{\beta}$
when $\alpha$ is a limit ordinal, and construct $p_{\alpha+1} 
\leq p_{\alpha}$ such that $p_{\alpha + 1} \cap \,^{\alpha}2 = 
p_{\alpha} \cap \,^{\alpha}2$.  For each $\alpha$, we satisfy the
following requirements:

{\bf (3.6)} \ For all $s \in p_{\alpha}$, if $o(s) < \alpha$ 
then:
\roster
\item"{(i)}"  If $s$ is a limit of splitting nodes in $p_{\alpha}$
		and $o(s) \in S_{\xi}$, then $s$ does not split in
		$p_{\alpha}$.
\item"{(ii)}"  If $s$ does not split in $p_{\alpha}$, then $(p_{\alpha})_s$
		decides the value of $\dof(o(s))$.
\item"{(iii)}"  If $s$ splits in $p_{\alpha}$, let $\gamma_s$ be the
		least $\gamma$ such that $(p_{\alpha})_s$ does not
		decide $\dof(\gamma)$.  Then $(p_{\alpha})_{s^{\frown}0}$
		and $(p_{\alpha})_{s^{\frown}1}$ decide $\dof(\gamma_s)$
		in opposite ways, and both $(p_{\alpha})_{s^{\frown}0}$
		and $(p_{\alpha})_{s^{\frown}1}$ have stems of length
		greater than $\gamma_s$.
\endroster

Note that if $p_{\alpha}$ satisfies (iii) for a given $s$, then 
every $p_{\beta}$, $\beta > \alpha$, satisfies (iii) for this $s$,
with the same $\gamma_s$.  Also (by induction on $o(s)$), we have
$\gamma_s \geq o(s)$.  Clearly, if $\alpha$ is a limit ordinal and
each $p_{\beta}$, $\beta < \alpha$, satisfies (3.6), then $p_{\alpha}$
also satisfies (3.6).  We show below how to obtain $p_{\alpha +1}$
when we have already constructed $p_{\alpha}$.

Now let $q = \bigcap_{\alpha < \kappa} p_{\alpha}$, and let us
verify that $q$ satisfies (3.5).  So let $\alpha$ be a limit ordinal,
and let $Z_{\alpha} = \{t \in q : t \text{ is a limit of splitting
nodes in  $q$ and }$ $o(t) = \alpha\}$.  If $t \in Z_{\alpha}$, then
$t$ is a limit of splitting nodes of $p_{\alpha}$.  It follows from
(3.6) (ii) and (iii) that $(p_{\alpha})_t$ decides $\dof \upharpoonright
\alpha$, and we let $h_t$ be this sequence.  If $t_1 \neq t_2$ are
in $Z_{\alpha}$, let $s = t_1 \cap t_2$.  By (3.6) (iii) we have
$\gamma_s < \alpha$ (because there exist $s_1$ and $s_2$ such that $s
\subset s_1 \subset t_1$, $s \subset s_2 \subset t_2$ and both $s_1$
and $s_2$ split in $p_{\alpha}$).  It follows that $h_{t_1}\neq
h_{t_2}$.  If $\alpha \in S_{\xi}$ and $s \in Z_{\alpha}$, then by
(3.6) (i), $s$ does not split in $p_{\alpha+1}$ and so $(p_{\alpha +
1})_s$ decides $\dof(\alpha)$; we let $i_s$ be this value.  These
$h_t$ and $i_s$ satisfy (3.5) for the condition $q$.

It remains to show how to obtain $p_{\alpha+ 1}$ from $p_{\alpha}$.
Thus assume that $p_{\alpha}$ satisfies (3.6).  First let $r \leq
p_{\alpha}$ be the following condition such that $r \cap \,^\alpha 2
= p_{\alpha} \cap \,^\alpha 2$: If $\alpha \notin S_{\xi}$ let
$r = p_{\alpha}$; if $\alpha \in S_{\xi}$, consider all $s \in p_{\alpha}$
with $o(s) = \alpha$ that are limits of splitting nodes, and replace
each $(p_{\alpha})_s$ by a stronger condition $r(s)$ such that $s$
does not split in $r(s)$.  For all other $s \in p_{\alpha}$ with
$o(s) = \alpha$, let $r(s) = (p_{\alpha})_s$.  Let $r$ be the amalgamation
of the $r(s)$; the tree $r$
is a condition because $s(\xi) = 1$ for all $s \in p_{\alpha}$
with $o(s) = \alpha$.

Now consider all $s \in r$ with $o(s) = \alpha$.  If $s$ does not
split in $r$, let $t$ be the successor of $s$ and let $q(t) \leq 
r_t$ be some condition that decides $\dof(\alpha)$.  If $s$ splits
in $r$, let $t_1$ and $t_2$ be the two successors of $s$, and let
$\gamma_s$ be the least $\gamma$ such that $\dof(\gamma)$ is not
decided by $r_s$.  Let $q(t_1) \leq r_{t_1}$ and $q(t_2) \leq 
r_{t_2}$ be conditions that decide $\dof(\gamma_s)$ in opposite
ways, and such that they have stems of length greater than 
$\gamma_s$.

Now we let $p_{\alpha + 1}$ be the amalgamation of all the
$q(t)$, $q(t_1)$, $q(t_2)$ into $r$.  Clearly, $p_{\alpha + 1}
\cap \,^\alpha 2 = r \cap \,^\alpha 2 = p_{\alpha} \cap \,^\alpha 2$.
The condition $p_{\alpha + 1}$ satisfies (3.6) (i) because $p_{\alpha}
\leq r$.  It satisfies (ii) because if $s$ does not split and
$o(s) = \alpha$, then $(p_{\alpha +1})_s = q(t)$ where $t$ is the
successor of $s$.  Finally, it satisfies (iii), because if $s$ splits
and $o(s) = \alpha$, then $(p_{\alpha + 1})_{s^{\frown}0} = q(t_1)$ and
$(p_{\alpha + 1})_{s^{\frown}1} = q(t_2)$ where $t_1$ and $t_2$ are the
two successors of $s$.  \pf

We now complete the proof that $B(P)$ is rigid.

\proclaim
{Theorem 3.7}  The complete Boolean algebra $B(P)$ has no
nontrivial automorphism.
\endproclaim

\demo
{Proof}  Assume that $\pi$ is a nontrivial automorphism of $B(P)$.
By Lemma 3.2 there exist conditions $p$ and $q$ with incomparable
stems $s$ and $t$ such that $\pi(p)$ and $q$ are compatible.
Let $t_0 = s \cap t$ and let $\xi = o(t_0)$.  Hence $t_0^{\frown} 0$
and $t_0^{\frown}1$ are stems of the two conditions and by Lemma 3.4,
one forces $\varphi(\xi)$ and the other forces $\neg \varphi(\xi)$.
This is a contradiction because $\pi(p)$ forces the same 
sentences that $p$ does, and $\pi(p)$ is compatible with $q$.\pf
\enddemo

\enddocument